\documentclass[11pt]{article}
\usepackage[latin1]{inputenc}
\usepackage{amsmath}
\usepackage{amsfonts}
\usepackage{amssymb,amsthm, stmaryrd}
\usepackage{graphicx}
\usepackage{xcolor}
\usepackage{multicol}
\usepackage{soul}
\usepackage[pagewise, displaymath, mathlines]{lineno}

\newcommand*\patchAmsMathEnvironmentForLineno[1]{%
  \expandafter\let\csname old#1\expandafter\endcsname\csname #1\endcsname
  \expandafter\let\csname oldend#1\expandafter\endcsname\csname end#1\endcsname
  \renewenvironment{#1}%
     {\linenomath\csname old#1\endcsname}%
     {\csname oldend#1\endcsname\endlinenomath}}%
\newcommand*\patchBothAmsMathEnvironmentsForLineno[1]{%
  \patchAmsMathEnvironmentForLineno{#1}%
  \patchAmsMathEnvironmentForLineno{#1*}}%
\AtBeginDocument{%
\patchBothAmsMathEnvironmentsForLineno{equation}%
\patchBothAmsMathEnvironmentsForLineno{align}%
\patchBothAmsMathEnvironmentsForLineno{flalign}%
\patchBothAmsMathEnvironmentsForLineno{alignat}%
\patchBothAmsMathEnvironmentsForLineno{gather}%
\patchBothAmsMathEnvironmentsForLineno{multline}%
}

\newtheorem{theorem}{Theorem}[section]

\newtheorem{proposition}[theorem]{Proposition}

\newtheorem{problem}{Problem}
\newtheorem{corollary}[theorem]{Corollary}

\newtheorem{remark}[theorem]{Remark}

\newcommand{\Z}{\mbox{${\mathbb Z}$}}
\newcommand{\N}{\mbox{${\mathbb N}$}}

\title{Erd\H{o}s-Ginzburg-Ziv type generalizations for linear equations and linear inequalities in three variables}

\begin{document}

\maketitle

\begin{center}

\begin{multicols}{2}

Mario Huicochea\\[1ex]
{\small Instituto de Matem\'aticas\\
UNAM Juriquilla\\
Quer\'etaro, Mexico\\
dym@cimat.mx}\\[2ex]

\columnbreak

Amanda Montejano\\[1ex]
{\small UMDI, Facultad de Ciencias\\
UNAM Juriquilla\\
Quer\'etaro, Mexico\\
amandamontejano@ciencias.unam.mx}\\[4ex]

\end{multicols}

\end{center}

\begin{abstract}
For any linear inequality in three variables $\mathcal{L}$, we determine (if it exist) the smallest integer $R(\mathcal{L}, \Z/3\Z)$ such that: for every mapping $\chi :[1,n] \to \{0,1,2\}$,  with $n\geq R(\mathcal{L}, \Z/3\Z)$, there is a solution $(x_1,x_2,x_3)\in [1,n]^3$ of $\mathcal{L}$ with $\chi(x_1)+\chi(x_2)+\chi(x_3)\equiv 0$ (mod $3$). Moreover,  we prove that  $R(\mathcal{L}, \Z/3\Z)=R(\mathcal{L}, 2)$, where  $R(\mathcal{L}, 2)$ denotes the classical $2$-color Rado number, that is, the  smallest integer (provided it exist) such that for every $2$-coloring of $[1,n]$, with $n\geq R(\mathcal{L}, 2)$, there exist  a monochromatic solution of $\mathcal{L}$. Thus, we get an Erd\H{o}s-Ginzburg-Ziv type generalization for all lineal inequalities in  three variables having  a solution in the positive integers. We also show a number of families of linear equations in  three variables  $\mathcal{L}$ such that they do not admit such Erd\H{o}s-Ginzburg-Ziv type generalization, named $R(\mathcal{L}, \Z/3\Z)\neq R(\mathcal{L}, 2)$. At the end of this paper some questions are proposed.
\end{abstract}

\section{Introduction}

In this paper we investigate colorings of sets of natural numbers. We denote by $[a,b]$ the interval of natural numbers $\{x\in \N\,:\,a\leq x\leq b\}$, and by $[a,b]^k$ the set of vectors $(x_1,x_2,...,x_k)$ where $x_i\in [a,b]$ for each $1\leq i\leq k$. An $r$-coloring of $[1,n]$ is a function $\chi:[1,n]\to [0,r-1]$. Given an $r$-coloring of $[1,n]$, a vector $(x_1,x_2,...,x_k)\in[1,n]^k$ is called \emph{monochromatic} if all its entries received the same color, \emph{rainbow} if all its entries received pairwise distinct colors, and \emph{zero-sum} if $\sum_{i=1}^k\chi(x_i)\equiv 0$ (mod $r$). 

For a Diophantine system of equalities (or inequalities) in $k$ variables  $\mathcal{L}$, we denote by $R(\mathcal{L}, r)$ the classical \emph{$r$-color Rado number}, that is, the  smallest integer, provided it exist, such that for every $r$-coloring of $[1,n]$, with $n\geq R(\mathcal{L}, r)$, there exist  $(x_1,x_2,...,x_k)\in [1,n]^k$ a solution of $\mathcal{L}$ which is monochromatic. Rado numbers have been widely studied for many years (see for instance \cite{LR}). When studying the existence of zero-sum solutions, it is common to refer to an $r$-coloring as a $(\Z/r\Z)$-coloring. In this setting,  Bialostocki, Bialostocki and Schaal \cite{bbs} started the study of the parameter $R(\mathcal{L}, \Z/r\Z)$ defined as the smallest integer, provided it exist, such that for every $(\Z/r\Z)$-coloring of $[1,R(\mathcal{L}, \Z/r\Z)]$ there exist a zero-sum solution of $\mathcal{L}$. Recently, Robertson and other authors studied the same parameter concerning different  equations or systems of equations, \cite{R18}, \cite{R20}, \cite{RRS}.

We shall note that, if $\mathcal{L}$ is a system of equalities (or inequalities)  in $k$ variables, then
\begin{equation}\label{eq:basic}
R(\mathcal{L}, 2)\leq R(\mathcal{L}, \Z/k\Z) \leq R(\mathcal{L}, k),
\end{equation} 
where the first inequality follows since, in particular, a $(\Z/k\Z)$-coloring that uses only colors $0$ and $1$ is a $2$-coloring where a zero-sum solution is a monochromatic solution; the second inequality of (\ref{eq:basic})  follows since any monochromatic solution of $\mathcal{L}$ in a $k$-coloring of $[1,R(\mathcal{L}, k)]$ is a zero-sum solution too.
In view of the Erd\H{o}s-Ginzburg-Ziv theorem \cite{EGZ}, the authors of \cite{bbs} state that a system $\mathcal{L}$ admits an EGZ-\emph{generalization} if $R(\mathcal{L}, 2)= R(\mathcal{L}, \Z/k\Z)$. For example, it is not hard to see that the system $AP(3): x+y=2z$ , $x<y$, admits an EGZ-generalization while  the Schur equation, $x+y=z$, does not. More precisely, we have that $$9=R(AP(3),2)=R(AP(3),\Z/3\Z)=9,$$ and   $$5=R(x+y=z,2)<R(x+y=z,\Z/3\Z)=10,$$ where  $R(AP(3),2)$ and  $R(x+y=z,2)$ are the well known van der Waerden number for $3$-term arithmetic progressions concerning  two colors and the Schur number concerning two colors respectively, while $R(x+y=z,\Z/3\Z)=10$ can be found in \cite{R18} and $R(AP(3),\Z/3\Z)=9$ can be found in \cite{R20}. 
In \cite{bbs} the authors consider the systems of  inequalities
\[\mathcal{L}_1: \sum_{i=1}^{k-1}x_i<x_k,\]
and
\[\mathcal{L}_2: \sum_{i=1}^{k-1}x_i<x_k\hspace{.2cm}, \hspace{.2cm}x_1<x_2< \cdots < x_k, \] 
proving that $\mathcal{L}_2$ admits an  EGZ-generalization for $k$ prime,   and $\mathcal{L}_1$ admits an  EGZ-generalization for any $k$, particularly, $R(\mathcal{L}_1, 2)= R(\mathcal{L}_1, \Z/k\Z)=k^2-k+1$ (see \cite{bbs}). In this paper we provide analogous results concerning any linear inequality on three variables. More precisely, let  $a,b,c,d\in \Z$, such that $abc\neq 0$. Then we consider,

\[\mathcal{L}_3: ax+by+cz+d<0.\]

We  prove that $\mathcal{L}_3$  admits an  EGZ-generalization for every set of integers $\{a,b,c,d\}$ such that the corresponding $2$-color Rado number exists. Moreover, we determine, in each case, such Rado numbers (see Theorem \ref{thm:L3}).

Note that, as we investigate linear systems, $\mathcal{L}$,  in three variables, to have an EGZ-generalization means that $$R(\mathcal{L},2)=R(\mathcal{L},\Z/3\Z),$$ and the parameter $R(\mathcal{L},\Z/3\Z)$ is defined as  the smallest integer, provided it exist, such that for every $f:[1,R(\mathcal{L}, \Z/3\Z)]\to \{0,1,2\}$ there exist a zero-sum (mod $3$) solution of $\mathcal{L}$ which, in this case, is either a monochromatic or a rainbow solution of $\mathcal{L}$. Therefore, the study of   $R(\mathcal{L},\Z/3\Z)$ is considered as a canonical Ramsey problem.


The paper is organized as follows. In Section \ref{sec:L3}, we find the explicit values of $R(\mathcal{L}_3, \Z/3\Z)$ and $R(\mathcal{L}_3, 2)$ (whenever $\mathcal{L}_3$ has solutions in the positive integers) in terms of the coefficients of $\mathcal{L}_3$. As a corollary, we get that $\mathcal{L}_3$ admits an EGZ-generalization in this case. In Section \ref{sec:negative} we provide some negative results, that is, we exhibit families of linear equations in three variables which admit no EGZ-generalization. In Section \ref{sec:existence}, we talk about $r$-regular linear equations and the families $\mathcal{F}_k$ and $\mathcal{F}_{\mathbb{Z}/k\mathbb{Z}}$. At the end of this section, we give some some problems related with these families.

\section{The $2$-color Rado numbers for $\mathcal{L}_3$}\label{sec:L3}
In this section we prove that any linear inequality on three variables, $\mathcal{L}_3$,  for which the corresponding $2$-color Rado number exists, admits an EGZ-generalization. We also determine the value of such Rado numbers depending on the coefficients of $\mathcal{L}_3$.

We will repeatedly use the following fact.

\begin{remark}\label{rem:AB}
Let $A$ and $B$ be integers such that $A<0$. Then, $$A\left(  \left\lfloor \frac{B}{-A} \right\rfloor +1\right)+B<0\leq A  \left\lfloor \frac{B}{-A} \right\rfloor +B.$$
\end{remark}

\begin{proof}
The first inequality follows since $\frac{B}{-A} - \left\lfloor \frac{B}{-A} \right\rfloor<1$, equivalently $\frac{B}{-A}< 1+ \left\lfloor \frac{B}{-A} \right\rfloor$ and,  multiplying both sides by $A$ (which is negative) we obtain $-B>A\left(  \left\lfloor \frac{B}{-A} \right\rfloor +1\right)$ from which it follows the claim. The second inequality holds true since $\left\lfloor \frac{B}{-A} \right\rfloor \leq \frac{B}{-A} $ and,  multiplying both sides by $A$ (which is negative) we obtain $A \left\lfloor \frac{B}{-A} \right\rfloor \geq A  \frac{B}{-A}=-B$ from which it follows the claim. 
\end{proof}

Recall that, for integers $a$, $b$, $c$ and $d$, such that $abc\neq 0$,
\[\mathcal{L}_3: ax+by+cz+d<0.\]

\begin{theorem}\label{thm:L3}
Let  $a,b,c,d\in \Z$, such that $abc\neq 0$, $a\leq b\leq c$, and define $\sigma=a+b+c+d$. If $\mathcal{L}_3$ has a solution in the positive integers, then

\[R(\mathcal{L}_3, 2)= R(\mathcal{L}_3, \Z/3\Z)= \left\{ \begin{array}{ll}
                  1 & \mbox{ if } \sigma<0,   \vspace{.3cm}   \\         \vspace{.3cm}  
                   \left\lfloor \frac{d}{-a-b-c} \right\rfloor  +1   & \mbox{ if } \sigma\geq 0  \mbox{ and } a\leq b\leq c<0,\\  \vspace{.3cm}  
  \left\lfloor \frac{c \left(\left\lfloor \frac{c+d}{-a-b} \right\rfloor  +1\right)+d}{-a-b} \right\rfloor  +1   & \mbox{ if } \sigma\geq 0  \mbox{ and } a\leq b<0<c,\\             \vspace{.3cm}       
   \left\lfloor \frac{\left(b+c\right)\left(\left\lfloor \frac{b+c+d}{-a} \right\rfloor  +1\right)+d}{-a} \right\rfloor  +1   & \mbox{ if } \sigma\geq 0  \mbox{ and } a<0<b\leq c,\\        
\end{array}\right.\]
\end{theorem}

\begin{proof}
Let $\{a,b,c,d\}$ be a set of integers such that $\mathcal{L}_3$ has some (integer) positive solution. If $a+b+c+d<0$ then $(1,1,1)$ is a monochromatic solution of $\mathcal{L}_3$ and so $R(\mathcal{L}_3, 2)= R(\mathcal{L}_3, \Z/3\Z)=1$. If $a+b+c+d\geq 0$, then necessarily some of the coefficients, $a$, $b$ or $c$, must be negative (otherwise, for all $x,y,z$ positive integers, $ax+by+cz+d\geq a+b+c+d\geq 0$ and $\mathcal{L}_3$ would have no solutions in the positive integers). Thus, assuming that $a+b+c+d\geq 0$, we consider three cases.\\

\noindent
\emph{Case 1. } Assume that $ a\leq b\leq c<0$. Define $k_0=\big\lfloor\frac{d}{-a-b-c}\big\rfloor+1$. First note that, since $a+b+c+d\geq 0$ and $-a-b-c>0$, then $k_0>1$. Observe now that, for any $x,y,z\in[1,k_0-1]$,
\begin{align*}
ax+by+cz+d&\geq a(k_0-1)+b(k_0-1)+c(k_0-1)+d\\
&=(a+b+c)\left\lfloor\frac{d}{-a-b-c}\right\rfloor+d\geq 0,
\end{align*}
where the last inequality follows by taking $A=a+b+c<0$ and $B=d$ in Remark~\ref{rem:AB}. Then, we conclude that $\mathcal{L}_3$ has no solution in $[1,k_0-1]$. On the other hand,   
\begin{align}\label{eq:3k_0}
ak_0+bk_0+ck_0+d&=(a+b+c)\left(\bigg\lfloor\frac{d}{-a-b-c}\bigg\rfloor+1\right)+d<0,
\end{align}
where the inequality follows by Remark \ref{rem:AB} (taking again $A=a+b+c<0$ and $B=d$).
From (\ref{eq:3k_0}), we conclude that $(k_0,k_0,k_0)$ is a  solution of $\mathcal{L}_3$, and so any coloring of $[1,k_0]$ will contain a monochromatic (zero-sum) solution of $\mathcal{L}_3$. Hence, $R(\mathcal{L}_3, 2)= R(\mathcal{L}_3, \Z/3\Z)=k_0$.\\

\noindent
\emph{Case 2. } Assume that $ a\leq b<0<c$. Define the function
\begin{equation*}
\psi:\Z\longrightarrow \Z,\qquad  \psi(x)=\bigg\lfloor\frac{cx+d}{-a-b}\bigg\rfloor+1,
\end{equation*}
and set $k_1=\psi(1)$ and $k_2=\psi(k_1)$. First note that, since $a+b+c+d\geq 0$ and $-a-b>0$ then $k_1>1$ and, as $\psi$ is a nondecreasing function, then $1<k_1\leq k_2$. From (\ref{eq:basic}), it suffices to show that
\begin{equation}
\label{lem3:L_3:E1}
k_2\leq R(\mathcal{L}_3, 2)
\end{equation}
 and 
 \begin{equation}
\label{lem3:L_3:E2}
R(\mathcal{L}_3, \Z/3\Z)\leq k_2. 
\end{equation}

To show (\ref{lem3:L_3:E1}), we exhibit a $2$-coloring of $[1,k_2-1]$ without  monochromatic solutions of $\mathcal{L}_3$.   Define $\chi_1:[1,k_2-1]\to\{0,1\}$ as 
\[\chi_1(x)=\left\{ \begin{array}{ll}
                  0 & \mbox{ if } 1\leq x\leq k_1-1,   \\ 
 		1 & \mbox{ if } k_1\leq x\leq k_2-1.  
\end{array}\right.\]
Note that, for all $x,y,z\in [1,k_1-1]$,
\begin{align}
\label{lem3:L_3:E3}
ax+by+cz+d&\geq a(k_1-1)+b(k_1-1)+c+d\nonumber\\
&=(a+b)\left\lfloor\frac{c+d}{-a-b}\right\rfloor+c+d\geq 0,
\end{align}
where the last  inequality follows by taking $A=a+b<0$ and $B=c+d$ in Remark~\ref{rem:AB}. Also, for all $x,y,z\in [k_1,k_2-1]$,
\begin{align}
\label{lem3:L_3:E4}
ax+by+cz+d&\geq a(k_2-1)+b(k_2-1)+ck_1+d\nonumber\\
&=(a+b)\left\lfloor\frac{ck_1+d}{-a-b}\right\rfloor +ck_1+d\geq 0,
\end{align}
where the last inequality follows by taking $A=a+b<0$ and $B=ck_1+d$ in Remark \ref{rem:AB}.
From (\ref{lem3:L_3:E3}) and (\ref{lem3:L_3:E4}), we conclude that there are no  monochromatic solutions of $\mathcal{L}_3$ with respect to $\chi_1$, which completes the  proof of (\ref{lem3:L_3:E1}).

Now we prove (\ref{lem3:L_3:E2}). Let  $\chi:[1,k_2]\to \{0,1,2\}$ be an arbitrary coloring, and assume that $\chi$ contains no zero-sum solutions of $\mathcal{L}_3$. We will use two times the first inequality of Remark \ref{rem:AB}. First take  $A=a+b<0$ and $B=c+d$ to obtain
\begin{align}\label{eq:2k_1,1}
ak_1+bk_1+c+d&=(a+b)\left(\left\lfloor\frac{c+d}{-a-b}\right\rfloor+1\right)+c+d<0.
\end{align}
Now, take $A=a+b<0$ and $B=ck_1+d$ to obtain
\begin{equation}\label{eq:2k_2k_1}
ak_2+bk_2+ck_1+d=(a+b)\left(\left\lfloor\frac{ck_1+d}{-a-b}\right\rfloor+1\right)+ck_1+d<0.
\end{equation}
By (\ref{eq:2k_1,1}) we know that $(k_1,k_1,1)$ is a solution of $\mathcal{L}_3$ which, by assumption, cannot be zero-sum. Suppose, without loss of generality,  that $\chi(1)=0$ and $\chi(k_1)=1$. Next, we prove that $\chi(k_2)$ cannot be $0$, $1$ or $2$.
\begin{itemize}
\item By (\ref{eq:2k_2k_1}), we know that  $(k_2,k_2,k_1)$ is a solution of $\mathcal{L}_3$, and so $ \chi(k_2)\neq \chi(k_1)=1$. 
\item Since $c>0$ and $k_1>1$ then $ak_2+bk_2+c+d\leq ak_2+bk_2+ck_1+d$, which together with (\ref{eq:2k_2k_1}) implies that $(k_2,k_2,1)$ is a solution of $\mathcal{L}_3$. Thus, $ \chi(k_2)\neq \chi(1)=0$. 
\item Since  $a<0$ and $k_2>k_1$ then $ak_2+bk_1+c+d\leq ak_1+bk_1+c+d$, which together with  (\ref{eq:2k_1,1}) implies that $(k_2,k_1,1)$ is a solution of $\mathcal{L}_3$. Thus, $ \chi(k_2)\neq 2$.
\end{itemize}
This contradiction implies the existence of a zero-sum solution in any $(\Z/3\Z)$-coloring of $[1,k_2]$, and we completed the proof of (\ref{lem3:L_3:E2}).\\

\noindent
\emph{Case 3. } Assume that $ a<0<b\leq c$. Define the function
\begin{equation*}
\phi:\Z\longrightarrow \Z,\qquad  \phi(x)=\bigg\lfloor\frac{(b+c)x+d}{-a}\bigg\rfloor+1,
\end{equation*}
and set $k_3=\phi(1)$ and $k_4=\phi(k_3)$. First note that, since $a+b+c+d\geq 0$, then $k_3>1$ and, as $\phi$ is a nondecreasing function, then $1<k_3\leq k_4$.
From (\ref{eq:basic}), it is enough to show that
\begin{equation}
\label{lem2:L_3:E1}
k_4\leq R(\mathcal{L}_3, 2)
\end{equation}
 and 
 \begin{equation}
\label{lem2:L_3:E2}
R(\mathcal{L}_3, \Z/3\Z)\leq k_4. 
\end{equation}

To show (\ref{lem2:L_3:E1}), we exhibit a $2$-coloring of $[1,k_4-1]$ without  monochromatic solutions of $\mathcal{L}_3$.   Define $\chi_2:[1,k_4-1]\to\{0,1\}$ as 
\[\chi_2(x)=\left\{ \begin{array}{ll}
                  0 & \mbox{ if } 1\leq x\leq k_3-1,   \\ 
 		1 & \mbox{ if } k_3\leq x\leq k_4-1.  
\end{array}\right.\]
Note that, for all $x,y,z\in [1,k_3-1]$,
\begin{align}
\label{lem2:L_3:E3}
ax+by+cz+d&\geq a(k_3-1)+b+c+d\nonumber\\
&=a\left\lfloor\frac{b+c+d}{-a}\right\rfloor+b+c+d\geq 0,
\end{align}
where the last inequality follows by taking $A=a<0$ and $B=b+c+d$ in Remark~\ref{rem:AB}. Also, for all $x,y,z\in [k_3,k_4-1]$,
\begin{align}
\label{lem2:L_3:E4}
ax+by+cz+d&\geq a(k_4-1)+bk_3+ck_3+d\nonumber\\
&=a\left\lfloor\frac{(b+c)k_3+d}{-a}\right\rfloor +bk_3+ck_3+d\geq 0,
\end{align}
where the last inequality follows by taking $A=a<0$ and $B=bk_3+ck_3+d$ in Remark~\ref{rem:AB}.
From (\ref{lem2:L_3:E3}) and (\ref{lem2:L_3:E4}), we conclude that there are no  monochromatic solutions of $\mathcal{L}_3$ with respect to $\chi_2$, which completes the  proof of (\ref{lem2:L_3:E1}).

Now we prove (\ref{lem2:L_3:E2}). Let  $\chi:[1,k_4]\to \{0,1,2\}$ be an arbitrary coloring, and assume that $\chi$  contains no zero-sum solutions of $\mathcal{L}_3$. 
We will use two times the first inequality of Remark \ref{rem:AB}. First take $A=a<0$ and $B=b+c+d$ to obtain
\begin{align}\label{eq:k_3,1,1}
ak_3+b+c+d&=a\left(\bigg\lfloor\frac{b+c+d}{-a}\bigg\rfloor+1\right)+b+c+d<0.
\end{align}
Now take $A=a<0$ and $B=(b+c)k_3+d$ to obtain 
\begin{equation}\label{eq:k_42k_3}
ak_4+bk_3+ck_3+d=a\left(\left\lfloor\frac{(b+c)k_3+d}{-a}\right\rfloor+1\right)+(b+c)k_3+d<0.
\end{equation}
By (\ref{eq:k_3,1,1}) we know that $(k_3,1,1)$ is a solution of $\mathcal{L}_3$ which, by assumption, cannot be zero-sum. Suppose, without loss of generality,  that $\chi(1)=0$ and $\chi(k_3)=1$. Next, we prove that $\chi(k_4)$ cannot be $0$, $1$ or $2$.
\begin{itemize}
\item By (\ref{eq:k_42k_3}), we know that  $(k_4,k_3,k_3)$ is a solution of $\mathcal{L}_3$, and so $ \chi(k_4)\neq \chi(k_3)=1$. 
\item  Since $c\geq b>0$ and $k_3>1$ then $ak_4+b+c+d\leq ak_4+bk_3+ck_3+d$, which together with (\ref{eq:k_42k_3}) implies that $(k_4,1,1)$ is a solution of $\mathcal{L}_3$. Thus, $ \chi(k_4)\neq \chi(1)=0$. 
\item Since  $c>0$ and $k_3>1$ then $ak_4+bk_3+c+d\leq ak_4+bk_3+ck_3+d$, which together with  (\ref{eq:k_42k_3}) implies that $(k_4,k_3,1)$ is a solution of $\mathcal{L}_3$. Thus, $ \chi(k_4)\neq 2$. 
\end{itemize}
This contradiction implies the existence of a zero-sum solution in any $(\Z/3\Z)$-coloring of $[1,k_4]$, and we completed the proof of (\ref{lem2:L_3:E2}).
\end{proof}

As an immediate consequence of Theorem \ref{thm:L3} we conclude the following.

\begin{corollary}\label{cor:L3}
Let  $a,b,c,d\in \Z$, such that $abc\neq 0$. If
$\mathcal{L}_3: ax+by+cz+d<0$
has a solution in the positive integers, then $\mathcal{L}_3$ admits an EGZ-generalization. 
\end{corollary}


\section{Negative results}\label{sec:negative}

In this section we exhibit different families  of linear equations in three variables which admit no EGZ-generalization. In other words, we study equations, $\mathcal{L}$,  where $R(\mathcal{L},2) \neq R(\mathcal{L},\Z/3\Z)$. Naturally, we focus our attention in equations such that both $R(\mathcal{L},2)$ and $R(\mathcal{L},\Z/3\Z)$ exist. Although Rado's Theorem characterizes the equations $\mathcal{L}$ such that $R(\mathcal{L},2)$ exists, there is a small number of families of equations where the value $R(\mathcal{L},2)$ is explicitly known, see \cite{libro}, \cite{sur}. In this section we develop some ideas to compare $R(\mathcal{L},2)$ and $R(\mathcal{L},\Z/3\Z)$ for some  equations, and then we get some applications to show that $\mathcal{L}$ does not admit EGZ-generalization. 

\begin{theorem}\label{thm:d}
Let $a,b,c,d$ be integers where $a,b, c$ are odd and $d$ is even, such that both  $R\left(ax+by+cz=\frac{d}{2}, 2\right)$ and $R(ax+by+cz=d, \Z/3\Z)$ exist. Then $$2R\left(ax+by+cz=\frac{d}{2}, 2\right)\leq R(ax+by+cz=d, \Z/3\Z).$$
\end{theorem}

\begin{proof}
Abbreviate writting $R:=R\left(ax+by+cz=\frac{d}{2}, 2\right)$. Let $\chi_0:[1,R-1]\rightarrow \{0,1\}$ be a coloring such that $ax+by+cz=\frac{d}{2}$ has not monochromatic solutions with respect to $\chi_0$. Define
\begin{equation*}
\chi:[1,2R-1]\longrightarrow \{0,1,2\},\qquad \chi(n)= \left\{ \begin{array}{ll}
                  \chi_0\left(\frac{n}{2}\right) & \mbox{ if $n$ is even}   \\         
                2  &    \mbox{ if $n$ is odd}    
\end{array}\right.
\end{equation*}
To prove the claim of the theorem, it is enough to show that $ax+by+cz=d$ has not zero-sum solutions with respect to $\chi$. Let $(x_0,y_0,z_0)$ be a solution of $ax+by+cz=d$. Since $d$ is even and $a,b,c$ are odd, we have that either the $3$ entries of $(x_0,y_0,z_0)$ are even or exactly one of the entries is even.

First assume that the $3$ entries of $(x,y,z)$ are even.  Then $\chi(x_0)=\chi_0\left(\frac{x_0}{2}\right)$, $\chi(y_0)=\chi_0\left(\frac{y_0}{2}\right)$ and $\chi(z_0)=\chi_0\left(\frac{z_0}{2}\right)$; since $ax+by+cz=\frac{d}{2}$ has not monochromatic solutions with respect to $\chi_0$, $\left(\frac{x_0}{2}, \frac{y_0}{2}, \frac{z_0}{2}\right)$ is not monochromatic. Therefore  $(x_0,y_0,z_0)$ is not a zero-sum solutions with respect to $\chi$. 
 
Now assume that exactly one of the entries of $(x,y,z)$ is even; without loss of generality assume that $x_0$ is even. Then $\chi(x_0)=\chi_0\left(\frac{x_0}{2}\right)$, $\chi(y_0)=2$ and $\chi(z_0)=2$. This means that $(x_0,y_0,z_0)$ is not a zero-sum solutions with respect to $\chi$.
\end{proof}

The next result is an immediate corollary of Theorem \ref{thm:d}.

\begin{corollary}\label{cor:odds}
Let  $\mathcal{L}$ be the equation $ax+by+cz=0$, and assume that $R(\mathcal{L}, \Z/3\Z)$ exist. Then, $\mathcal{L}$ admits no EGZ-generalization if $a$, $b$ and $c$ are odd integers.
\end{corollary}

\begin{proof}
By  taking  $d=0$ in Theorem \ref{thm:d}, we   conclude that $2R(\mathcal{L}, 2)\leq R(\mathcal{L}, \Z/3\Z)$, and so $R(\mathcal{L}, 2)\neq R(\mathcal{L}, \Z/3\Z)$.  
\end{proof}

Also Theorem \ref{thm:d} provides some applications for non-homogeneous linear equations.
\begin{corollary}\label{cor:nh1}
Let  $d$ be a negative even integer. Then, the equation  $x+y-z=d$ admits no EGZ-generalization.
\end{corollary}
\begin{proof}
From \cite[Thm. 9.14]{libro}, we have that 
\begin{equation*}
R(x+y-z=d,2)=5-4d, 
\end{equation*}
 and 
 \begin{equation*}
R\left(x+y-z=\frac{d}{2},2\right)=5-2d. 
\end{equation*}
On the other hand, Theorem \ref{thm:d} leads to 
\begin{equation*}
2R\left(x+y-z=\frac{d}{2},2\right)\leq R(x+y-z=d,\Z/3\Z). 
\end{equation*}
Hence
\begin{align*}
R(x+y-z=d,\Z/3\Z)&\geq 2R\left(x+y-z=\frac{d}{2},2\right)\\
&=2(5-2d)\\
&>5-4d\\
&=R(x+y-z=d,2).\qedhere
\end{align*}
\end{proof}

\begin{corollary}\label{cor:nh2}
Let  $d$ be a positive integer congruent to $6,8$ or $0$ modulo $10$. Then, the equation  $x+y-z=d$ admits no EGZ-generalization.
\end{corollary}
\begin{proof}
On the one hand, we have from \cite[Thm. 9.15]{libro} that 
\begin{equation*}
R(x+y-z=d,2)=d-\left\lceil\frac{d}{5}\right\rceil+1, 
\end{equation*}
 and 
 \begin{equation*}
R\left(x+y-z=\frac{d}{2},2\right)=\frac{d}{2}-\left\lceil\frac{\frac{d}{2}}{5}\right\rceil+1. 
\end{equation*}
On the other hand, Theorem \ref{thm:d} leads to 
\begin{equation*}
2R\left(x+y-z=\frac{d}{2},2\right)\leq R(x+y-z=d,\Z/3\Z). 
\end{equation*}
Thus, since $d$ is congruent to $6,8$ or $0$ modulo $10$, we get that
\begin{align*}
R(x+y-z=d,\Z/3\Z)&\geq 2R\left(x+y-z=\frac{d}{2},2\right)\\
&=2\left(\frac{d}{2}-\left\lceil\frac{\frac{d}{2}}{5}\right\rceil+1\right)\\
&>d-\left\lceil\frac{d}{5}\right\rceil+1\\
&=R(x+y-z=d,2).\qedhere
\end{align*}
\end{proof}
\section{Other directions}\label{sec:existence}
A linear homogenous equation is called \emph{$r$-regular} if every $r$-coloring of $\N$ contains a monochromatic solution of it (equivalently, an equation $\mathcal{L}$ is called  $r$-regular if $R(\mathcal{L},r)$ exist). A linear homogenous equation is called \emph{regular} if it is $r$-regular for all positive integers $r$. 
Denote by $\mathcal{F}_r$ the family of linear homogenous equations which are $r$-regular. For equations on $k\geq 3$ variables, Rado completely determined $\mathcal{F}_2$:  it is the set of  equations, $\sum_{i=1}^k c_ix_i=0$ for which there exist  $i,j\in\{1,\cdots,k\}$ such that  $c_i<0$  and  $c_j>0$ (see, for instance \cite{libro}). For other values of $r\in \Z^+$, the family $\mathcal{F}_r$ is not characterized. Rado's Single Equation Theorem states that a  linear homogenous equation on $k\geq 2$ variables, $\sum_{i=1}^k c_ix_i=0$  ($c_i$'s are non-zero integers),  is regular if and only if there exist a non-empty $D\subseteq \{1,...,k\}$ such that  $\sum_{d\in D} c_d=0$.  
Naturally, $\mathcal{F}_{r+1} \subseteq \mathcal{F}_r$ for all $r\in \Z^+$. In his Ph.D. dissertation, Rado conjectured that, for all $r\in \Z^+$, there are equations that are $r$-regular but not $(r+1)$-regular. This conjecture was solved by Alexeev and Tsimerman in 2010 \cite{alexev}, where they confirm that $\mathcal{F}_{r+1} \subsetneq \mathcal{F}_r$.  For any $k\in\N$, define $\mathcal{F}_{\mathbb{Z}/k\mathbb{Z}}$ to be the family of linear homogeneous  equations in $k$ variables, $\mathcal{L}$, for which $R(\mathcal{L},\Z/k\Z)$  exist. By (\ref{eq:basic}) we know that $$\mathcal{F}_3\subseteq \mathcal{F}_{\mathbb{Z}/3\mathbb{Z}}\subseteq \mathcal{F}_2.$$ We will show that 
\begin{equation}
\label{eq:nonin}
\mathcal{F}_3\subsetneq \mathcal{F}_{\mathbb{Z}/3\mathbb{Z}}.
\end{equation}
For all $n\in\N$,  denote by  $\mathrm{ord}_2(n)$ the maximum $m\in\mathbb{Z}$  such that  $2^m$  divides $n$. First note that the equation $x+2y-4z=0$ is not in $\mathcal{F}_3$ since the coloring
 \[\chi:\N\to\{0,1,2\},\qquad \chi(n)=\left\{ \begin{array}{ll}
                  0 & \mbox{ if } \mathrm{ord}_2(n)\equiv 0 \mod 3,   \\ 
                   1 & \mbox{ if } \mathrm{ord}_2(n)\equiv 1 \mod 3,   \\ 
 	 2 & \mbox{ if } \mathrm{ord}_2(n)\equiv 2 \mod 3.
\end{array}\right.\]
has not monochromatic solutions of it. The next proposition implies that $x+2y-4z=0$ is in $\mathcal{F}_{\mathbb{Z}/3\mathbb{Z}}$ and therefore (\ref{eq:nonin}) holds.
\begin{proposition}
The equation $x+2y-4z=0$ satisfies $R(x+2y-4z=0, \mathbb{Z}/3\mathbb{Z})=8$.
\end{proposition}

\begin{proof}
Let  $\mathcal{L}$ be the equation $x+2y-4z=0$. First we show that $R(\mathcal{L}, \mathbb{Z}/3\mathbb{Z})\geq 8$. Let $\chi:[1,8]\rightarrow \{0,1,2\}$ be a coloring. We assume that there are not zero-sum solutions of $\mathcal{L}$ with respect to $\chi$, and we will get a contradiction.  Assume without loss of generality that $\chi(1)=0$. Since $(2,1,1)$ is a solution of $\mathcal{L}$, $\chi(2)\neq 0$; assume without loss of generality that $\chi(2)=1$. Note that $(2,3,2)$,  $(4,2,2)$ and $(4,4,3)$ are solutions  of $\mathcal{L}$. Thus either $\chi(3)=0$ and $\chi(4)=2$, or $\chi(3)=2$ and $\chi(4)=0$.  Notice that $(6,3,3)$ and   $(4,6,4)$ are solutions of $\mathcal{L}$ so $\chi(6)=1$. Since $(6,5,4)$ and   $(2,5,3)$ are solutions of $\mathcal{L}$, we get $\chi(5)=1$. For any value of $\chi(8)$, we obtain a zero-sum solution inasmuch as $(8,2,3)$, $(8,4,4)$ and $(8,6,5)$ are solutions of $\mathcal{L}$ and this is the desired contradiction.

On the other hand, $R(\mathcal{L}, \mathbb{Z}/3\mathbb{Z})>7$ since there are not  zero-sum solutions with respect to the coloring
\begin{equation*}
\chi:[1,7]\longrightarrow \{0,1,2\},\qquad \chi(n)= \left\{ \begin{array}{ll}
                  0 & \mbox{ if $n\in\{1,4,7\}$}   \\         
                1  &    \mbox{ if $n\in \{2,5,6\}$}  \\
                 2  &    \mbox{ if $n=3$},   
\end{array}\right.
\end{equation*}
 and this completes the proof.
\end{proof}
From the previous discussion, we know that $\mathcal{F}_3\subsetneq \mathcal{F}_{\mathbb{Z}/3\mathbb{Z}}\subseteq \mathcal{F}_2$. A natural question arises from this chain.
\begin{problem}
\label{pro:1}
Is it true that $ \mathcal{F}_{\mathbb{Z}/3\mathbb{Z}}\subsetneq \mathcal{F}_2$?
\end{problem} 
From (\ref{eq:basic}) we know that $\mathcal{F}_k\subsetneq \mathcal{F}_{\mathbb{Z}/k\mathbb{Z}}$ for all $k\geq 3$. Thus it would be interesting to know if there are $k\in\N$ such that the equality is achieved.
\begin{problem}
\label{pro:2}
 For all $k\geq 3$,  $\mathcal{F}_k \subsetneq \mathcal{F}_{\mathbb{Z}/k\mathbb{Z}}$?
\end{problem}
Finally we know  that $\mathcal{F}_k\subsetneq \mathcal{F}_{\mathbb{Z}/k\mathbb{Z}}$ and $\mathcal{F}_{k}\subseteq \mathcal{F}_{k-1}$ for all $k\geq 3$. However we do not know whether there is a relation between $\mathcal{F}_{k-1}$ and $\mathcal{F}_{\mathbb{Z}/k\mathbb{Z}}$.
\begin{problem}
\label{pro:3}
For all $k\geq 3$, $\mathcal{F}_{\mathbb{Z}/k\mathbb{Z}}\subseteq \mathcal{F}_{k-1} $?
\end{problem}

\section{Acknowledgements}

The present work was initiated during the workshop Zero-Sum Ramsey Theory: Graphs, Sequences
and More (19w5132). We thank the facilities provided by the Banff International Research Station ``Casa Matem\'{a}tica Oaxaca".
 This research was partially supported by  PAPIIT project IN116519.




\end{document}